\documentclass{amsart}
\usepackage{amsmath, amsfonts, amssymb, latexsym}

\newtheorem{theorem}{Theorem}

\newtheorem{conjecture}[theorem]{Conjecture}
\newtheorem{corollary}[theorem]{Corollary}

\newtheorem{lemma}[theorem]{Lemma}
\newtheorem{notation}[theorem]{Notation}

\newtheorem{proposition}[theorem]{Proposition}
\newtheorem{remark}[theorem]{Remark}

\numberwithin{equation}{section}
\newcommand{\be}[1]{\Rand{\vspace{0,6cm}\tt #1}\begin{equation}\label{#1}}
\begin{document}
\title[Branching Brownian motion with self-interaction]
      {The center of mass for spatial branching processes and
       an application for self-interaction}
\author{J\'anos Engl\"ander}
\address{Department of Statistics and Applied Probability\\
         University of California, Santa Barbara\\
         Santa Barbara 93106}
\email{englander@pstat.ucsb.edu}
\urladdr{http://www.pstat.ucsb.edu/englander/}
\thanks{Research supported by  NSA  grant
~H982300810021} 
\date{8/21/08}

\begin{abstract}
In this paper we prove that the center of mass of a supercritical
branching-Brownian motion, or that of a supercritical super-Brownian
motion tends to a limiting position almost surely, which, in a sense
complements a result of Tribe on the final behavior of a critical
super-Brownian motion. This is shown to be true also  for a model
where branching Brownian motion is modified by attraction/repulsion
between particles.

We then put this observation together with the description of the
interacting system as viewed from its center of mass, and get the
following asymptotic behavior: the system asymptotically becomes a
branching Ornstein Uhlenbeck process (inward for attraction and
outward for repulsion), but
\begin{enumerate}
\item  the origin is shifted to a random point which has normal distribution, and
\item the Ornstein Uhlenbeck particles are not independent but constitute a system with a degree of freedom which is less
by their number by precisely one.
\end{enumerate}
\end{abstract}

\maketitle

%\tableofcontents
\section{Introduction}\label{S:intro}

We start with some basic notation.
\begin{notation}\rm In this paper $\mathcal{M}_f(\mathbb R^d)$ and
$\mathcal{M}_1(\mathbb R^d)$ denote the space of finite measures and
the space of probability measures, respectively, on $\mathbb{R}^d$.
For $\mu\in\mathcal{M}_f(\mathbb R^d)$, we define
$\|\mu\|:=\mu(\mathbb R^d)$.
\end{notation}
\subsection{A model with self-interaction}\label{details}
Consider a dyadic (i.e. precisely two offspring replaces the parent)
branching Brownian motion (BBM) in $\mathbb R^d$ with unit time
branching and with the following interaction between particles: if
$Z$ denotes the process and $Z_t^i$ is the $i$th particle, then
$Z_t^i$ feels the drift
$$\frac{1}{n_t}\sum_{1\le j \le n_t}\gamma \left(Z_t^j-\cdot\right),$$ where $\gamma\neq 0$ ,
that is the particle's infinitesimal generator is
\begin{equation}\label{generator}
\frac{1}{2}\Delta +\frac{1}{n_t}\sum_{1\le j \le n_t}\gamma
\left(Z_t^j-x\right)\cdot \nabla.
\end{equation}  (Here and in the sequel, $n_t$ is a shorthand for $2^{\lfloor t\rfloor}$, where
$\lfloor t\rfloor$ is the integer part of $t$.) If $\gamma>0$, then
this means \emph{attraction}, if $\gamma<0$, then it means
\emph{repulsion}.

To be a bit more precise, we can define the process by induction as
follows. $Z_0$ is a single particle at the origin. In the time
interval $[m,m+1)$ we define a system of $2^m$ \emph{interacting
diffusions}, by the following system of SDE's:
$$Z_{\tau}^i=W_{\tau}^i+\frac{\gamma}{2^m}\int_0^{\tau}\sum_{1\le j\le 2^m}(Z_s^i-Z_s^i)\, \mathrm{d}s,\
\tau\in[0,1);\ i=1,2,\dots, 2^m,$$ where $\{W_{\tau}^i\,
,\tau\in[0,1)\}$ are independent Brownian motions for $i=1,2,\dots,
2^m$ and they start at the position of their parents at the end of
the previous step (at time $m-0$). Existence and uniqueness follows
from the fact that these $2^m$ interacting diffusions can be
considered as a single $2^m d$-dimensional Brownian motion with
drift $\mathbf{b}:\mathbb{R}^{2^m d}\rightarrow \mathbb{R}^{2^m d}$
where
\begin{eqnarray*} \mathbf{b}\big(x_1,x_2,...,x_d,x_{1+d},x_{2+d},...,x_{2d},...,x_{1+(2^m-1)d},x_{2+(2^m-1)d},...,x_{2^md})\\
=:(\beta_1,\beta_2,...,\beta_{2^md})^T
\end{eqnarray*}
satisfies
\begin{eqnarray*}
&&\beta_1=2^{-m}(x_1+x_{1+d}+...+x_{1+(2^m-1)d})-x_1,\\
&&\beta_2=2^{-m}(x_2+x_{2+d}+...+x_{2+(2^m-1)d})-x_2,\\
&&...\\
&&\beta_d=2^{-m}(x_d+x_{2d}+...+x_{2^md})-x_d,\\
&&\beta_{1+d}=2^{-m}(x_1+x_{1+d}+...+x_{1+(2^m-1)d})-x_{1+d},\\
&&\beta_{2+d}=2^{-m}(x_2+x_{2+d}+...+x_{2+(2^m-1)d})-x_{2+d},\\
&&...\\
&&\beta_{2d}=2^{-m}(x_d+x_{2d}+...+x_{2^md})-x_{2d},\\
&&...\\
&&\beta_{1+(2^m-1)d}=2^{-m}(x_1+x_{1+d}+...+x_{1+(2^m-1)d})-x_{1+(2^m-1)d},\\
&&\beta_{2+(2^m-1)d}=2^{-m}(x_2+x_{2+d}+...+x_{2+(2^m-1)d})-x_{2+(2^m-1)d},\\
&&...\\
&&\beta_{2^m d}=2^{-m}(x_d+x_{2d}+...+x_{2^md})-x_{2^m d}.
\end{eqnarray*}
Existence and uniqueness is guaranteed by the Lipschitzness
(linearity) of $b$.

\begin{remark}\rm
It may seem natural to replace the interaction we defined by the
gravitational force between particles\footnote{I.e. the forces are
given by Newton's law of universal gravitation: they vary as the
inverse square of the distance between the particles.}, however this
would lead to a randomized version of the (notoriously difficult)
`n-body problem.'$\hfill\diamond$
\end{remark}
\subsection{Results on the self-interacting model}
We are interested in the long time behavior of $Z$, and also,whether
we can say something about the number of particles in a given
compact set for $n$ large. In the sequel we will use the standard
notation $\langle Z_t,g\rangle:=\sum_{i=1}^{n_t} g(Z_t^{i}).$

In this paper we will show (Theorem \ref{asymptotic.behavior}) that
$Z$ asymptotically becomes a branching Ornstein Uhlenbeck process
(inward for attraction and outward for repulsion), but
\begin{enumerate}
\item  the origin is shifted to a random point which has $d$-dimensional normal
distribution $\mathcal{N}^d(0,2)$, and
\item the Ornstein Uhlenbeck particles are not
independent but constitute a system with a degree of freedom which
is less by their number by precisely one.
\end{enumerate} For the
local behavior we formulate and motivate a conjecture (Conjecture
\ref{br.diff}).
\subsection{An extension of Tribe's result on critical super-Brownian motion}
In the proof of Theorem \ref{asymptotic.behavior} we will first show
that $\overline{Z}_t:=\frac{1}{n_t}\sum_{i=1}^{n_t}Z_t^i$,  the
center of mass for $Z$ satisfies
$\lim_{t\to\infty}\overline{Z}_t=N$, where $N\sim
\mathcal{N}^d(0,2)$. In fact, the proof will reveal that
$\overline{Z}$ moves like a Brownian motion, which is nevertheless
slowed down tending to a final limiting location (see Lemma
\ref{cofm} and its proof).

Since this is also true for $\gamma=0$ (BBM with unit time branching
and no self-interaction), our first natural question is whether we
can prove a similar result for the supercritical super-Brownian
motion.

Another motivation for the same goal is as follows. Tribe \cite{T92}
proved that a \emph{critical} super-Brownian motion near its
extinction time $\xi$ behaves like a single Brownian path stopped at
$\xi$. More precisely, $X_t\rightarrow \delta_F$ as $t\to\infty$
a.s. in the weak topology, where $F$ is a $d$-dimensional random
variable and its distribution is the same as that of a Brownian
motion at time $\xi$.

We would like to extend Tribe's result to the supercritical
super-Brownian motion $X$. We are interested in whether we can
obtain a similar result on the \emph{survival set}. Of course, then
$X_t$ does not shrink to a point in any sense, however, we may hope
to get an analogous result regarding the center of mass, defined as
$$\overline{X}_t:=\frac{1}{\|X_t\|}\int_{\mathbb R^d}x\, X_t(\mathrm{d}x)=
\langle x/\|X_t\|,X_t\rangle.$$ (Since $f(x)=x$ is not a bounded
function we may not hope to use Tribe's techniques though.)

In the sequel we will prove that $\overline{X}$ will  be a
\emph{time changed Brownian motion on a finite} (but random)
\emph{time interval}.

At the beginning of this subsection we referred to a result on
$\overline{Z}$. In fact we now see that the results on
$\overline{X}$ and $\overline{Z}$ are completely analogous: in both
cases the center of mass is a  Brownian motion slowed down in such a
way that the time interval $[0,\infty)$ is compressed into a finite
one. A slight difference is that, in case of $Z$, the terminal time
is deterministic, because so is the offspring distribution. (The
terminal time is $t=2$.)

Let $X$ be the $(\frac{1}{2}\Delta,\beta,\alpha; \mathbb
R^d)$-superdiffusion with $\alpha,\beta>0$ (supercritical
super-Brownian motion). Let $P$ denote the corresponding
probability. Let us restrict $\Omega$ to the  survival set
$$S:=\{\omega\in\Omega\mid X_t(\omega)>0,\ \forall t>0\}.$$ Since $\beta>0$, $P_{\mu}(S)>0$ for all $\mu\neq \mathbf{0}$.

Our main result is that the center of mass  for $X$ stabilizes as
$t\to\infty$, and furthermore the path of the center of mass is a
finite piece of a Brownian path (with a different time
parametrization).
\begin{theorem}\label{sbm}
Let $\alpha,\beta>0$ and let $\overline{X}$ denote the center of
mass process  for the $(\frac{1}{2},\beta,\alpha; \mathbb
R^d)$-superdiffusion $X$.  Then, on $S$,
\begin{itemize}
\item[(i)]
$\overline{X}_t$ converges $P_{\delta_{x}}$-almost surely as
$t\to\infty$.
\item[(ii)] In fact, the finite  path  $\{\overline{X}_t\}_{t\ge 0}$
is the same as the path of a time changed Brownian motion on
$[0,T_{\infty})$, where $T_{\infty}$ is a positive and finite random
variable. More precisely, there exists a $d$-dimensional Brownian
motion $B$ (on an enlarged space) such that $\overline{X}=B\circ T$
where $T:\mathbb{R}_+\to\mathbb{R}_+$ is a random time change
satisfying $\lim_{t\to\infty} T(t)<\infty\ P_{\delta_{x}}-a.s.$
\end{itemize}
\end{theorem}
\begin{remark}\rm

$ $
\begin{itemize}

\item [(a)]  A heuristic argument for (i) is as follows. Obviously, the center of mass is invariant under $H$-transforms whenever $H$ is
spatially (but not temporarily) constant. Let $H(t):=e^{-\beta t}.$
Then $X^H$ is a $(\frac{1}{2}\Delta,0,e^{-\beta t}\alpha; \mathbb
R^d)$-superdiffusion, that is, a critical super-Brownian motion with
a clock that is slowing down. Therefore, heuristically it seems
plausible that $\overline{X^H}$, the center of mass for the
transformed process stabilizes, because this is obviously true in
case of extinction, and otherwise the center of mass under the heat
flow does not move.

\item[(b)] From (ii) one can easily conclude some of the \emph{a.s. path properties} of $\overline{X}$.
For example, since the $p$-variation ($p>0$) is invariant under
changing the parametrization,  we know  that any segment of the path
has infinite total variation, but the whole path has finite
quadratic variation a.s. $\hfill \diamond$
\end{itemize}
\end{remark}
\section{The mass center stabilizes}
Notice that \begin{equation}\label{notice} \frac{1}{n_t}\sum_{1\le
j\le n_t} \left(Z_t^j-Z_t^i\right)=\frac{1}{n_t}\left(\sum_{1\le j
\le n_t} Z_t^j-n_t Z_t^i\right)=\overline{Z}_t-Z_t^i,
\end{equation}
and so \emph{the net attraction pulls the particle towards the
center of mass} (net repulsion pushes it away from the center of
mass). Thus the following lemma is relevant:
\begin{lemma}[Mass center stabilizes]\label{cofm}
Let $\overline{Z}_t:=\frac{1}{n_t}\sum_{i=1}^{n_t}Z_t^i$, that is,
$\overline{Z}$ is the  center of mass for $Z$. Then
$\lim_{t\to\infty}\overline{Z}_t=N$, where $N\sim
\mathcal{N}^d(0,1)$.
\end{lemma}
\begin{proof}
For $t\in [m,m+1)$, there are $2^m$ particles moving around. If
$t,t+\Delta t    \in [m,m+1)$, then $Z_t^i$ moves as a Brownian
motion plus a vector
$$\gamma (\overline{Z}_t-Z_t^i)\Delta t,$$
so $$Z^i_{t+\Delta t}=Z_t^i+\gamma (\overline{Z}_t-Z_t^i)\Delta t+
B_0^i(\Delta t),$$ where the $\{B_0^i(s),\ s\ge 0;\ i=1,2,...,2^m\}$
are independent Brownian motions starting at the origin. Thus
$$\overline{Z}_{t+\Delta t}=2^{-m}\sum_{i=1}^{2^{m}}Z^i_{t+\Delta t}=
 \overline{Z}_t+\gamma 2^{-m}\left(2^m\overline{Z}_t-\sum_{i=1}^{2^{m}}
Z_t^i\right)\Delta t+2^{-m}\bigoplus_{i=1}^{2^{m}}B_0^i(\Delta t),$$
where $\bigoplus$ denotes independent sum. Since
$2^m\overline{Z}_t-\sum_{i=1}^{2^{m}} Z_t^i=0$, we obtain that for
$0\le \tau<1$,
\begin{equation}\label{com.moves}
\overline{Z}_{m+\tau}=\overline{Z}_m +2^{-m}\bigoplus_{i=1}^{2^{m}}
B_0^i(\tau).
\end{equation}
Let
$\widehat{B}^{(m)}(\tau):=2^{-m/2}\bigoplus_{i=1}^{2^{m}}B_0^i(\tau)$.
Using induction, we obtain that
$$\overline{Z}_t=\widehat B^{(0)}(1)\oplus\frac{1}{\sqrt{2}}\widehat B^{(1)}(1)\oplus\dots\frac{1}{2^{k/2}}\widehat
B^{(k)}(1)\dots\oplus\frac{1}{\sqrt{2^{\lfloor t\rfloor-1}}}\hat
B^{(2^{\lfloor t\rfloor-1})}(1)\oplus\frac{1}{\sqrt{n_{t}}}\widehat
B^{(n_t)}(1),$$ (it is easy to check that, as the notation suggests,
the summands are independent) which, using Kolmogorov's Three Series
Theorem converges almost surely. (We will denote the a.s. limit by
$N$.) On the other hand, since $\widehat{B}^{(m)}$ is a Brownian
motion, we can  apply Brownian scaling and get
$2^{-m/2}\widehat{B}^{(m)}(\tau)\stackrel{d}{=}
W^{(m)}\left(\frac{\tau}{2^m}\right),$ where $W^{(m)},\ m\ge 1$ are
independent Brownian motions. We have
$$\overline{Z}_{m+\tau}\stackrel{d}{=}
\overline{Z}_m\oplus W^{(m)}\left(\frac{\tau}{2^m}\right),$$ and so,
\begin{eqnarray*}
\widehat B^{(0)}(1)\oplus\frac{1}{\sqrt{2}} \widehat
B^{(1)}(1)\oplus\dots\oplus\frac{1} {\sqrt{2^{\lfloor
t\rfloor-1}}}\widehat B^{(2^{\lfloor
t\rfloor-1})}(1)\oplus\frac{1}{\sqrt{n_{t}}}\widehat
B^{(n_t)}(1)\\
\stackrel{d}{=}W^{(1)}(1)\oplus
W^{(2)}\left(\frac{1}{2}\right)\oplus\dots\oplus  W^{(2^{\lfloor
t\rfloor-1})}\left(\frac{1}{2^{\lfloor t\rfloor-1}}\right)\oplus
W^{(n_t)}\left(\frac{\tau}{n_t}\right),
 \end{eqnarray*}
where $\tau:=t-n_t$. Since the summands are independent on the right
hand side,  $N$ has the same distribution as a Brownian motion at
$t=2$, that is, $N\sim\mathcal{N}(0,2)$.
\end{proof}
For another proof see the remark after Lemma \ref{normal}.
\begin{remark}\rm
It is interesting to note that $\overline{Z}$ \emph{is in fact a
Markov process}. Indeed, the distribution of $\overline{Z}_t$
($t=m+\tau,\ 0\le\tau<1$)  conditional on $\mathcal{F}_m$ is the
same as conditional on  $Z_m$, because $Z$ itself is a Markov
process. But the distribution of $\overline{Z}_t$ only depends  on
${Z}_m$ through $\overline{Z}_m$, as
$$\overline{Z}_t\stackrel{d}{=}\overline{Z}_m\oplus W^{(2^{m})}\left(\frac{\tau}{2^{m}}\right),$$ whatever ${Z}_m$ is.$\hfill\diamond$
\end{remark}
We will also need the following fact later.
\begin{lemma} The coordinate processes of $Z$ are independent.\label{indep}
\end{lemma}
We leave the simple proof to the reader.

\section{Normality via decomposition}
We will need the following result. The \emph{decomposition}
appearing in the proof will also be useful. Recall that $m=\lfloor
t\rfloor$.
\begin{lemma}\label{normal} $(Z_t^1,Z_t^2,...,Z_t^{2^m})$ is joint normal for all $t\ge 0$.
\end{lemma}

\noindent{\bf Proof of lemma.} By Lemma \ref{indep}, we may assume
that $d=1$. We prove the statement by induction.

For $m=1$ it is trivial.

Suppose that the statement is true for $m-1$. Since
$(Z_{m-1}^1,Z_{m-1}^2,...,Z_{m-1}^{2^m})$ is normal, we can consider
it just as well as a $2^m$ dimensional degenerate normal at the
instant of the fission of the particles. Indeed, the vector
$$(Z_{m-1}^1,Z_{m-1}^1,Z_{m-1}^2,Z_{m-1}^2,...,Z_{m-1}^{2^m},Z_{m-1}^{2^m})$$
has the same distribution on the $2^{m-1}$ dimensional subspace
$$S:=\{x\in \mathbb{R}^{2^{m}}\mid
x_1=x_2,x_3=x_4,...,x_{2^{m}-1}=x_{2^{m}}\}$$ of
$\mathbb{R}^{2^{m}}$ as the vector
$\sqrt{2}(Z_{m-1}^1,Z_{m-1}^2,...,Z_{m-1}^{2^m})$ on
$\mathbb{R}^{2^{m-1}}$. (The reader can easily visualize this for
$m=2$: the distribution of $(Z_1^1,Z_1^1)$ is clearly $\sqrt{2}$
times the distribution of a Brownian particle at time $1$, i.e.
$\mathcal{N}(0,\sqrt{2})$ on the line $x_1=x_2$.)

Since the convolution of normals is normal, therefore, by the Markov
branching property, it is enough to prove the statement when the
$2^m$ particles start at the origin and the clock is reset:
$t\in[0,1).$

Define the $2^m$ dimensional process $Z^*$ on the time interval
$t\in[0,1)$ by
$$Z_t^*:=(Z_t^1,Z_t^2,...,Z_t^{2^m}),$$
starting at the origin. Because of the interaction between the
particles attracts the particles towards the center of mass, $Z^*$
is a Brownian motion with drift
$$\gamma \left[(\overline{Z}_t,\overline{Z}_t,...,\overline{Z}_t)-(Z_t^1,Z_t^2,...,Z_t^{2^m})\right].$$
Notice that this drift is orthogonal to the vector\footnote{For
simplicity, we use row vectors in this proof.} ${\mathbf
v}:=(1,1,...,1)$, that is, the vector
$(\overline{Z}_t,\overline{Z}_t,...,\overline{Z}_t)$ is nothing but
the orthogonal projection of $(Z_t^1,Z_t^2,...,Z_t^{2^m})$ to the
line of ${\mathbf v}$. This observation immediately leads to the
following decomposition. The process $Z^*$ can be decomposed into
two components:
\begin{itemize}
\item the component in the direction of ${\mathbf v}$ is a  Brownian motion
\item in the ortho-complement of ${\mathbf v}$, it is an independent Ornstein-Uhlenbeck process with parameter
$\gamma$.
\end{itemize}
It follows from this decomposition that $Z^*$ is Gaussian.
$\hfill\qed$
\begin{remark}\rm
Consider the Brownian component in the decomposition appearing in
the proof. Since, on the other hand,  this coordinate is
$2^{m/2}\overline{Z}_t$, using Brownian scaling, one obtains another
way of seeing that $\overline{Z}_t$ stabilizes at a position which
is distributed as the time $1+2^{-1}+2^{-2}+...+2^{-m}+...=2$ value
of a Brownian motion. (The decomposition shows this for $d=1$ and
then it is immediately upgraded to general $d$ by
independence.)$\hfill\diamond$
\end{remark}
\begin{corollary}[Asymptotics for finite subsystem]
Let $k\ge 1$ and consider the subsystem $(Z_t^1,Z_t^2,...,Z_t^k),\ t\ge m_0$ for  $m_0:= \lfloor\log k\rfloor+1$. (This
means that at time $m_0$ we pick $k$ particles and at every fission replace the parent particle by randomly picking one
of its two descendants.) Let the real numbers $c_1,...,c_k$  satisfy
\begin{equation}\label{two.cond}
\sum_{i=1}^k c_i=0,\ \sum_{i=1}^k c_i^2=1.
\end{equation} Define $\Psi^{(c_1,...,c_k)}_t:=\sum_{i=1}^k c_i Z_t^i$ and note that $\Psi_t$
is  invariant under the translations of the coordinate system. Let $\mathcal{L}_t$ denote its law.

For every $k\ge 1$ and $c_1,...,c_k$ satisfying (\ref{two.cond}), $\Psi^{(c_1,...,c_k)}$ is the same $d$-dimensional
Ornstein Uhlenbeck process corresponding to the operator $1/2\Delta -\gamma\nabla\cdot x$, and in particular,
$$\lim_{t\to\infty}\mathcal{L}_t=\mathcal{N}\left(0,\frac{d}{2\gamma}\right).$$
\end{corollary}
For example, taking $c_1=1/\sqrt{2},c_2=-1/\sqrt{2}$, we obtain that  when viewed from a tagged particle's position,
any given other particle moves as $\sqrt{2}$ times the above Ornstein Uhlenbeck process.

\noindent{\bf Proof.} By independence (Lemma \ref{indep}) it is enough to consider $d=1$. For $m$ fixed, consider the
decomposition appearing in the proof of Lemma \ref{normal} and recall the notation. By (\ref{two.cond}), whatever $m\ge
m_0$ is, the $2^m$ dimensional unit vector $$(c_1,c_2,...,c_k,0,0,...,0)$$ is orthogonal to the $2^m$ dimensional
vector ${\mathbf v}$. This means that $\Psi^{(c_1,...,c_k)}$ is a one dimensional projection of the Ornstein Uhlenbeck
component of $Z^*$, and thus it is itself a one dimensional Ornstein Uhlenbeck process (with parameter $\gamma$) on the
unit time interval.

Now, although as $m$ grows, the Ornstein Uhlenbeck components of $Z^*$ are defined on larger and larger spaces
($S\subset\mathbb{R}^{2^{m}}$ is a $2^{m-1}$ dimensional linear subspace), the projection onto the direction of
$(c_1,c_2,...,c_k,0,0,...,0)$ is always the same one dimensional Ornstein Uhlenbeck process, i.e. the different unit
time `pieces' of $\Psi^{(c_1,...,c_k)}$ obtained by those projections may be concatenated. $\hfill\qed$
\section{The interacting system as viewed from the center of mass}
Recall that by (\ref{com.moves}) the interaction has no effect on
the motion of $\overline{Z}$. Let us see now how the interacting
system looks like when viewed from $\overline{Z}$.
\subsection{The description of a single particle}\label{single}
Using our usual notation, assume that $t\in[m,m+1)$ and let
$\tau:=t-\lfloor t \rfloor$. Also recall that $\oplus$ denotes
independent sum. When viewed from $\overline{Z}$, the
 time $\tau$ relocation\footnote{I.e. the relocation between time
  $m$ and time $t$.} of a particle can be described as follows:
$$\Delta (Z^{1}_t-\overline{Z}_t)=\Delta Z^{1}_t-\Delta\overline{Z}_t=
B_0^{1}(\tau)-2^{-m}\bigoplus_{i=1}^{2^m} B_0^{i}(\tau)-\gamma(Z^{1}_t-\overline{Z}_t)\tau.$$ So if
$Y^{1}:=Z^1-\overline{Z}$, then
$$\Delta Y^{1}_t=
B_0^{1}(\tau)-2^{-m}\bigoplus_{i=1}^{2^m} B_0^{i}(\tau)-\gamma
Y^{1}_t\tau.$$ Clearly,
$$B_0^{1}(\tau)-2^{-m}\bigoplus_{i=1}^{2^m}
B_0^{i}(\tau)=\bigoplus_{i=2}^{2^m}
2^{-m}B_0^{i}(\tau)\oplus(1-2^{-m})B_0^1(\tau),$$ and thus the right
hand side is a Brownian motion with mean zero and variance
$d(1-2^{-m})\tau$. That is, denoting $\sigma_m:=1-2^{-m}$,
$$\Delta Y^{1}_t=
d\sigma_m W^1(\tau)-\gamma Y^{1}_t\tau,$$ where $W^1$ is a standard
Brownian motion.

We have thus obtained that (for fixed $m$) $Y^1$ corresponds to the
operator\footnote{Here $\Delta$ is the Laplace operator unlike in
previous lines where it denoted difference.}
$$\frac{1}{2}\sigma_m\Delta-\gamma x\cdot \nabla,$$ which is an
Ornstein-Uhlenbeck process on the time interval $[m,m+1)$. Since for
$m$ large $\sigma_m$ is close to one, the relocation viewed from the
center of mass is \emph{asymptotically governed by an O-U process
corresponding to $\frac{1}{2}\Delta-\gamma x\cdot \nabla.$}
\begin{remark}[Asymptotically vanishing correlation between driving
BM's]\rm Let $W^i=(W^{i,k},k=1,2,...,d)$. For $1\le i\neq j\le 2^m$,
we have
 \begin{eqnarray*}
&&E\left[\sigma_m W^{i,k}(\tau)\cdot\sigma_m W^{j,k}(\tau) \right]=\\
&&E\left[\left(B_0^{i,k}(\tau)-2^{-m}\bigoplus_{r=1}^{2^m}
B_0^{r}(\tau)\right)\left(B_0^{j,k}(\tau)-2^{-m}\bigoplus_{r=1}^{2^m}
B_0^{r}(\tau)\right)\right]=\\
&&2^{-m}\left[\text{Var}\left(B_0^{i,k}(\tau)\right)+
\text{Var}\left(B_0^{j,k}(\tau)\right)\right]-2^{-2m}\cdot
2^m\tau=(2^{1-m}-2^{-m})\tau=2^{-m}\tau,
 \end{eqnarray*}
that is, for $i\neq j$,
\begin{equation}\label{covariance}
E\left[W^{i,k}(\tau) W^{j,\ell}(\tau)
\right]=\delta_{k\ell}\cdot\frac{2^{-m}}{(1-2^{-m})^2}\,\tau.
\end{equation}
Hence \emph{the pairwise correlation decays to zero} as $t\to\infty$
 (recall that $m=\lfloor t\rfloor$ and $\tau=t-m\in [0,1)$).

And of course, for the variances we have
\begin{equation}\label{variance}
E\left[W^{i,k}(\tau) W^{i,\ell}(\tau) \right]=\delta_{k\ell}\cdot
(1-2^{-m})\,\tau,\ \text{for}\ 1\le i\le 2^m.\hspace{1cm}\diamond
\end{equation}
\end{remark}
\subsection{The description of the system; the
`degree of freedom'}  Fix $m\ge 1$ and for $t\in[m,m+1)$ let
$Y_t:=(Y_t^1,...,Y_t^{2^{m}})^T$, where $()^T$ denotes transposed.
(This is a vector of length $2^m$ where each component itself is a
$d$ dimensional vector; one can actually look at it as a $2^m\times
d$ matrix too.) We then have
$$\Delta Y_t=\sigma_m \left(W^1(\tau),...,W^{2^{m}}(\tau)\right)^T-\gamma Y_t\cdot \tau,$$
where $$W^i(\tau):=B_0^{i}(\tau)-2^{-m}\bigoplus_{j=1}^{2^m}
B_0^{j}(\tau),\ i=1,2,...,2^m$$  are mean zero Brownian motions with
correlation structure given by (\ref{covariance})-(\ref{variance}).

Just like at the end of subsection \ref{details}, we can consider
$Y$ as a single $2^md$-dimensional diffusion.  Each of its
components is an Ornstein-Uhlenbeck process with asymptotically unit
diffusion coefficient.

By independence, it is enough to consider the $d=1$ case, and so
from now on,  in this subsection we assume that $d=1$.

Let us first describe the distribution of $W_t$ for $t\ge 0$ fixed.
Recall that $\{B_0^i(s),\ s\ge 0;\ i=1,2,...,2^m\}$ are independent
Brownian motions starting at the origin. By definition, $W_t$ is a
$2^m$-dimensional multivariate normal:
\[ W_t=\left( \begin{array}{ccc}
1-2^{-m} & -2^{-m} & ... \  -2^{-m} \\
-2^{-m} & 1-2^{-m} & ... \  -2^{-m} \\
.\\
.\\
.\\
 -2^{-m} & -2^{-m} & ...  \  1-2^{-m} \end{array} \right)B_0(t)=:\mathbf{A}^{(m)} B_0(t)\]
However, since we are viewing the system from the center of mass, it
is a \emph{singular} multivariate normal. Its `true' dimension is
the rank of the matrix $\mathbf{A}^{(m)}$.
\begin{lemma}\label{rank}
rank$(\mathbf{A}^{(m)})=2^m-1.$
\end{lemma}

\noindent {\bf Proof.} We will simply write $\mathbf{A}$ instead of
$\mathbf{A}^{(m)}$. Since the columns of $\mathbf{A}$ add up to
zero, the matrix $\mathbf{A}$ is not of full rank: $r(\mathbf{A})\le
2^m-1.$ On the other hand,
$$2^m\mathbf{A}+ \left(
\begin{array}{ccc}
1 & 1 & ... \  1 \\
1 & 1 & ... \  1 \\
.\\
.\\
.\\
1 & 1 & ...  \ 1 \end{array} \right)=2^m \mathbf{I},
$$
where $\mathbf{I}$ is the $2^m$-dimensional unit matrix, and so by
subadditivity,
$$r(\mathbf{A})+1=r(2^m\mathbf{A})+1\ge 2^m.\hspace{1cm}\qed
$$
By Lemma \ref{rank}, $W_t$  is concentrated on the
$(2^m-1)$-dimensional linear subspace given by the orthogonal
complement of the vector $(1,1,...,1)^T$; in this $2^m-1$
dimensional subspace $W_t$ has \emph{non-singular} multivariate
normal distribution.  What this means is that even though
$W_1,W_2,...,W_{2^m}$ are not independent, their `degree of freedom'
is $2^m-1$, i. e. the  $2^m$-dimensional vector $W$ is
\emph{determined by $2^m-1$ independent components} (corresponding
to $2^m-1$ principal axes).

%Concerning the joint distribution at different times, let $k\ge 1$ and consider the $(k\times 2^m)$ random matrix
%$\mathcal{W}$ whose columns are the $2^m$-dimensional random (normal) vectors $W_{t_{1}},W_{t_{2}},\dots,W_{t_{k}}$. We
%have $$W_{t_{i}}=\mathbf{A}^{(m)} B_0(t_{i}),\ 1\le i\le k$$ and it is easy to show that the  random matrix
%$\mathcal{W}$ is normal, i.e. the joint distribution of its entries is $k 2^m$-dimensional normal. Indeed, ...

\section{Asymptotic behavior}
How can we put together that $\overline{Z}_t$ tends to a random
final position a.s. with the description of the system `as viewed
from $\overline{Z}_t$?'
\begin{lemma}\label{indep}
For $t\ge 0$, the random vector $Y_t$ is independent of  the path
$\{\overline{Z}_s\}_{s\ge t}$.
\end{lemma}

\noindent{\bf Proof.} First, for any $t>0$, $Y_t$ is independent of
$\overline{Z}_t$, because (assuming $d=1$) the vector
$$(\overline{Z}_t,
Z^1_t-\overline{Z}_t,Z^2_t-\overline{Z}_t,\dots,Z^{2^{m}}_t-\overline{Z}_t)^T$$ is normal (since it is a linear
transformation of the vector $(Z^1_t,Z^2_t,\dots,Z^{2^{m}}_t)^T$, which is normal by Lemma \ref{normal}), and so it is
sufficient to show that $\overline{Z}_t$ and $Z^i_t-\overline{Z}_t$ are uncorrelated for $1\le i\le 2^m$. But this is
obvious, because the random variables $Z^1_t,Z^2_t,...,Z^{2^{m}}_t$ are exchangeable and thus, denoting $n=2^m$,
\begin{eqnarray*}
E[\overline{Z}_t(Z^1_t-\overline{Z}_t)]=
E\left[\frac{Z^1_t+Z^2_t+...+Z^n_t}{n}\left(Z^1_t-\frac{Z^1_t+Z^2_t+...+Z^n_t}{n}\right)\right] = \\
 = \frac{1}{n}E\left[(Z^1_t)^2\right]+\frac{n-1}{n}E(Z^1_t
Z^2_t)-\frac{n}{n^2}E\left[(Z^1_t)^2\right]-\frac{n(n-1)}{n^2}E(Z^1_t
Z^2_t)=0.
\end{eqnarray*}
To complete the proof of the lemma, recall Lemma \ref{cofm} and its
proof and notice that the distribution of $\{\overline{Z}_s\}_{s\ge
t}$ only depends on its starting point $\overline{Z}_t$, as it is
that of a Brownian path appropriately slowed down, whatever $Y_t$
(or even $Z_t$) is. Since, as we have seen, $Y_t$ is independent of
$\overline{Z}_t$, we are done.$\hfill\qed$

\medskip
We have the following result on the asymptotic behavior of the
system.
\begin{theorem}\label{asymptotic.behavior}
The asymptotic behavior of $Z$ is that of a branching
Ornstein-Uhlenbeck process, but with the origin shifted by a random,
normally distributed vector. More precisely, for almost all $x$,
conditionally on $\lim_{t\to\infty}\overline{Z}_t=x$, one has
$Z_t=\overline{Z}_t\oplus Y_t$, where $Y$ is the branching
Ornstein-Uhlenbeck process of section \ref{single}.
\end{theorem}

\noindent{\bf Proof.} We already know from Lemma \ref{cofm} that the
almost sure limit $N:=\lim_{t\to\infty}\overline{Z}_t$ exists. Let
$$P^{x}(\cdot):=P(\cdot\mid N=x).$$ By Lemma \ref{indep}, $Y_t$ is
independent of $N$, and thus $P^{x}(Y_t\in \cdot)=P(Y_t\in \cdot)$
for almost all $x$. By definition, $Z_t=\overline{Z}_t+Y_t,$ and by
the above discussion, we can in fact write $Z_t=\overline{Z}_t\oplus
Y_t$ under $P^{x}$. Putting this together with the description of
$Y_t$ in section \ref{single}, we are done. $\hfill\qed$

We close this section with a conjecture on the local behavior of the
system.
\begin{conjecture}\label{br.diff}
Let $g\in C_c^+(\mathbb R^d)$.
\begin{itemize}
\item[(i)] If $\gamma>0$ (attraction), then there exists a
random variable $N\sim \mathcal{N}^d(0,2)$, such that, conditional
on $N=x_0$,
$$\lim_{n\to\infty}2^{-n}\langle Z_n,g\rangle=
\left\langle \left(\frac{\gamma}{\pi}\right)^{d/2} \exp\left(-\gamma
|\cdot-x_0|^2\right),g\right\rangle\ \ \mathrm{a.s.}$$ Also, there
exists a random variable $M\sim
\mathcal{N}^d\left(0,\left(2+\frac{1}{4\gamma^2}\right)\mathbf{I}_d\right)$
(with corresponding law $\mathbb{P}$), such that
$$\lim_{n\to\infty}2^{-n}E\langle Z_n,g\rangle=\mathbb{E}\left\langle
M,g\right \rangle\ \ \mathrm{a.s.}$$ (Here $\mathbf{I}_d$ denotes
the $d$-dimensional unit matrix.)
\item[(ii)]   If $\gamma<0$ (repulsion), then
$$\lim_{n\to\infty}2^{-n}\langle Z_n,g\rangle=
\langle 1,g\rangle\ \ \mathrm{a.s.}$$
\end{itemize}
\end{conjecture}
\begin{remark}\rm
The function $2+\frac{1}{4\gamma^2}$ in Conjecture \ref{br.diff}(i)
is monotone decreasing. The intuitive meaning is that stronger
attraction results in smaller  variance of the limiting
distribution.$\hfill\diamond$
\end{remark}

\noindent\textbf{Explanation of the conjecture:} Once we have
Theorem \ref{asymptotic.behavior}, we can try to put it together
with the Strong Law of Large Numbers for the local mass from
\cite{EHK08} for the process $Y$. So $N$ is the final position of
the center of mass.

First, let $\gamma>0$. If the components of $Y$ were independent and
the branching rate were exponential, Theorem 6  of \cite{EHK08}
would be readily applicable. However, since the $2^m$ components of
$Y$ are not independent (as we have seen, their degree of freedom is
$2^m-1$) and since, unlike in \cite{EHK08}, we now have unit time
branching, the method of \cite{EHK08} must be adapted to our
setting. This adaption would require some extra work.

Once this is done, however, it immediately follows that there exists
a random variable $N\sim \mathcal{N}^d(0,2)$, such that, conditional
on $N=x_0$,
$$\lim_{n\to\infty}2^{-n}\langle Z_n,g\rangle=
\left\langle \left(\frac{\gamma}{\pi}\right)^{d/2} \exp\left(-\gamma
|\cdot-x_0|^2\right),g\right\rangle\ \ \mathrm{a.s.}$$ Thus the
distribution of $\lim_{n\to\infty}2^{-n}\langle Z_n,g\rangle$ is the
convolution of $\mathcal{N}(0,2I_d)$ and
$\mathcal{N}(0,(2\gamma)^{-2}I_d)$, yielding the second statement in
the first part.

The $\gamma>0$ case is similar but since the limiting density is
translation invariant, i.e. Lebesgue (see Example 11 in
\cite{EHK08}), the final position of the center of mass plays no
role.
\section{Proof of Theorem \ref{sbm}}\label{Pf2}
(i) Since $\alpha,\beta$ are constant, the branching is independent
of the motion, and therefore $N$ defined by
$$N_t:={e^{-\beta t}\|X_t\|}$$ is a nonnegative martingale (positive on $S$) tending to a limit almost surely.
It is straightforward to check that it is uniformly bounded in $L^2$
and is therefore uniformly integrable (UI). Write
$$\overline{X}_t=\frac{e^{-\beta t}\langle
x,X_t\rangle}{e^{-\beta t}\|X_t\|}=\frac{e^{-\beta t}\langle
x,X_t\rangle}{N_t}.$$ We now claim that  $N_{\infty}>0$ a.s. on $S$.
Let $A:=\{N_{\infty}=0\}$. Clearly $S\complement\subset A$, and so
if we show that $P(A)=P(S\complement)$, then we are done. As is well
known, $P(S\complement)=e^{-\beta/\alpha}$. On the other hand, a
standard martingale argument (see the argument after formula (20) in
\cite{E08}) shows that $0\le u(x):=-\log P_{\delta_x}(A)$ must solve
the equation $$\frac{1}{2}\Delta u+\beta u-\alpha u^2=0,$$ but since
$P_{\delta_x}(A)=P(A)$ constant, therefore $-\log P_{\delta_x}(A)$
solves $\beta u-\alpha u^2=0$. Since $N$ is UI, no mass is lost in
the limit, giving $P(A)<1$. So $u>0$, which in turn implies that
$-\log P_{\delta_x}(A)=\beta/\alpha$.

Once we know that $N_{\infty}>0$ a.s. on $S$ , it is enough to focus
on the term $e^{-\beta t}\langle x,X_t\rangle$. Let $H(t):=e^{-\beta
t}$. Then $X^{H}$ is a $(\frac{1}{2}\Delta,0,e^{-\beta t}\alpha;
\mathbb R^d)$-superdiffusion, that is, a critical super-Brownian
motion with a clock that is slowing down. One can write $$e^{-\beta
t}\langle x,X_t\rangle=\langle x,X^{H}_t\rangle.$$ Define
$T_s={\mathcal S}^H_s:=
   e^{-\beta s}{\mathcal S}_s$;
then the semigroup $\{T_s\}_{s\ge 0}$ corresponds to Brownian
motion. In particular then
\begin{equation}\label{invariant} T_s [
\text{id}]=\text{id},
\end{equation}
where $\text{id}(x)=x.$ Therefore $\langle x, X^H_t\rangle$ is a
martingale.\footnote{It does not matter that the function is
unbounded and changes sign.} If we show that the martingale is UI,
we are done. It is enough to show that it is uniformly bounded in
$L^2$. To achieve this, define $g_n$ by $g_n(x)=|x|\cdot
\mathbf{1}_{\{|x|<n\}}.$ Then we have
$$E_{\delta_{0}}\langle x, X^H_t\rangle^2=E_{\delta_{0}}|\langle x, X^H_t\rangle |^2\le
E_{\delta_{0}}\langle |x|, X^H_t\rangle ^2,$$ and by the monotone
convergence theorem we can continue with
$$=\lim_{n\to\infty}E_{\delta_{0}}\langle g_n(x), X^H_t\rangle ^2.$$
Since $g_n$ is compactly supported, there is no problem to use the
moment formula and continue with
$$=\lim_{n\to\infty}\left( 1+\int_0^t \,\text{d}s\, e^{-\beta s} \langle \delta_{0}, T_s[\alpha
g_n^2]\rangle\right)=1+\alpha \lim_{n\to\infty}\int_0^t
\,\text{d}s\, e^{-\beta s} T_s[ g_n^2](0).$$ Recall that $\{T_s;
s\ge 0\}$ is the Brownian semigroup, that is,
$T_s[f](x)=\mathbb{E}_x (W_s)$, where $W$ is Brownian motion. Since
$g_n(x)\le |x|$, therefore we can trivially upper estimate the last
expression by
$$1+\alpha \lim_{n\to\infty}\int_0^t
\,\text{d}s\, e^{-\beta s} \mathbb{E}_0 (W^2_s)=1+\alpha \int_0^t
\,\text{d}s\, s e^{-\beta s}=1+\frac{\beta}{\alpha}.$$ Since this
upper estimate is independent of $t$, we are done:
$$\sup_{t\ge 0} E_{\delta_{0}}\langle x, X^H_t\rangle^2\le
1+\frac{\beta}{\alpha}.$$

(ii) Keeping in mind that  we are on the survival set
$S:=\{\omega\in\Omega\mid X_t(\omega)>0,\ \forall t>0\}$, we first
claim that $\overline{X}$ is an $(S,P,\{\mathcal{F}_t\}_{t\ge
0})$-martingale, where $\mathcal{F}_t:=\sigma(\{X_t;\ t\ge 0\})$.

Let us start with showing that $E|\overline{X}_t|<\infty$. To see
this, first recall that
$$\langle x,X_t\rangle/\|X_t\|=\frac{\langle x,X^H_t\rangle}{N_t}.$$
Thus,
$$E\left[|\overline{X}_t|;S\right]=E\left[\,N_t^{-1}|\langle
x,X^H_t\rangle|;S\right],$$ and since we have seen in part (i) that
$\langle x,X^H_t\rangle$ is (uniformly) bounded in $L^2$, by
Cauchy-Schwartz, it is enough to see that
$$E\left[N_t^{-2};S\right]=E[\|X^H_t\| ^{-2};S]<\infty,\ \forall
t>0,$$ i.e. that $E[\|X_t\| ^{-2};S]<\infty,\ \forall t>0$. This is
true because $\|X\|$ is a one-dimensional diffusion on $[0,\infty)$
with generator $x(\alpha \frac{d^2}{dx^2}+\beta\frac{d}{dx})$,
which, on $S$, tends to infinity, and therefore, by Fatou's Lemma,
$\lim_{t\to\infty}E[\|X_t\| ^{-2};S]=0.$ Hence $E[\|X_t\|
^{-2};S]<\infty$ for large $t$'s, but then by continuity, $E[\|X_t\|
^{-2};S]<\infty,\ \forall t>0$.

Next, since $E|\overline{X}_t|<\infty$, in fact $E\overline{X}_t=0$
for $t\ge 0$, because $\overline{X}_t$ is symmetrically distributed
($\overline{X}_t\stackrel{d}{=}-\overline{X}_t$). Indeed,
$$P\left(\frac{X_t}{\|X_t\|}\in B\right)=P\left(\frac{X_t}{\|X_t\|}\in B^*\right),
\ \forall B\in \mathcal{M}_1(\mathbb R^d),$$ where $B^*:=\{\mu*\mid
\mu\in B\}$ and $\mu^*(\cdot):=\mu(-\cdot)$.

We now show that $E(\overline{X}_t\mid
\mathcal{F}_s)=\overline{X}_s$ for $0\le s<t$. By the Markov
branching property,
$$E(\overline{X}_t\mid
\mathcal{F}_s)=E \frac{1}{\|X_s\|}\int _{\mathbb
R^d}\overline{X}_{t-s}^{\delta_{x}}\, X_s(\text{d}x),$$ where
$X^{\delta_{x}}$ are independent copies of super-Brownian motions
starting at $\delta_{x}$. Since $\overline{X}$ is mean zero,
shifting by $x$ and using Fubini's Theorem, we can continue the last
displayed formula with
$$=\frac{1}{\|X_s\|}\int _{\mathbb R^d}x\, X_s(\text{d}x)=\overline{X}_s.$$
Next, we show that $\overline{X}$ has continuous paths.

Let $$u(x,t):=x \mathbf{1}_{\{|x|\le
(\sqrt{2\beta}+\epsilon)t\}}+(\sqrt{2\beta}+\epsilon)t\mathbf{1}_{\{|x|>
(\sqrt{2\beta}+\epsilon)t\}},$$ where $\epsilon>0$ is fixed for the
rest of the argument.  By \cite{Ky04} it follows that there exists
an  a.s. finite random time $T=T(\omega)$ such that for all $t>T$,
$$\langle x,X_t\rangle=\langle u(\cdot, t),X_t\rangle.$$
(This is obviously true on the extinction set.) Hence, for $t>T$,
\begin{eqnarray*}\langle x,X_t\rangle-\langle x,X_{t+\Delta t}\rangle=\langle
u(\cdot,t),X_{t}\rangle-\langle u(\cdot,t+\Delta t),X_{t+\Delta
t}\rangle=  &&\\ \big(\langle u(\cdot,t),X_t\rangle-\langle
u(\cdot,t),X_{t+\Delta t}\rangle\big)+\big(\langle
u(\cdot,t),X_{t+\Delta t}\rangle-\langle u(\cdot,t+\Delta
t),X_{t+\Delta t}\rangle \big)&&\\=:I+II.&&
\end{eqnarray*}
Now $I\to 0$ as $\Delta t\downarrow 0$ a.s. by the continuity of the
paths (the continuity is in the weak topology and $u(\cdot, t)$ is a
bounded continuous function for all $t>0$), and $II\to 0$ as $\Delta
t\downarrow 0$ a.s., because
$$II=\langle \gamma(\Delta t),X_{t+\Delta t}\rangle,$$ where
$|\gamma(\Delta t)|\le (\sqrt{2\beta}+\epsilon)\Delta t\
\text{a.s}$. So
$$II\le (\sqrt{2\beta}+\epsilon)\Delta t\cdot \|X_{t+\Delta t}\|\ \text{a.s.},$$
and we are done since $$\lim_{\Delta t\downarrow 0}\|X_{t+\Delta
t}\|=\|X_{t}\|\ \text{a.s.},$$ again because of the path continuity.

This proves right continuity of the paths of $\overline{X}$; the
left continuity is similar.

Once we know that the center of mass is a continuous martingale, we
can utilize the multidimensional Dambis-Dubins-Schwarz
Representation Theorem, which we cite below. (For the theorem and
its proof see e.g. Theorem 18.4 in \cite{K02}.)

A $d$-dimensional continuous local $\mathcal{F}$-martingale
$M=(M^1,M^2,...,M^d)$ is said to be \emph{isotropic} if $[ M^1]=[
M^2]=...=[ M^d]$ holds almost surely, where $[ M^i]$ denotes the
total variation process for $M^i$.
\begin{proposition}[Theorem DDS]
Let $M$ be an isotropic continuous local $\mathcal{F}$-martingale in
$\mathbb{R}^d$ with $M_0$, and define
$$T_s:=\inf\{t\ge 0; [ M^1]_t>s\},\ \mathcal{G}_s:=
\mathcal{F}_{T_{s}},\ s\ge 0.$$ Then there exists in $\mathbb{R}^d$
a Brownian motion $B$ (with respect to a standard extension of
$\mathcal{G}$), such that a.s. $B=M\circ T$ on $[0,\langle
M^1\rangle_{\infty})$ and $M=B\circ[ M^1]$.
\end{proposition}
In order to finish the proof of Theorem \ref{sbm}(ii), first let
$d=1$. Since isotropy automatically holds,  if $M:=\langle x,
X^H\rangle$ and $T_s:=\inf\{t\ge 0; [ M]_t>s\}$, then by Theorem
DDS, there exists a one-dimensional Brownian motion $B$ (on an
enlarged space) such that $\overline{Z}=B\circ T$.

Now let $d\ge 2$. Note, that in our case
$$\langle x, X^H_t\rangle=\left(\langle x_1, X^H_t\rangle,\langle x_2,
X^H_t\rangle,...,\langle x_d, X^H_t\rangle\right),$$ if
$x=(x_1,x_2,...,x_d),$ and thus $M:=\langle x, X^H\rangle$ is
isotropic, because $\langle x_i, X^H\rangle$ is a time-changed
Brownian motion ($1\le i\le d$) by the already proven
one-dimensional case and so $[\langle x_i, X^H\rangle]_t=T_t$ for
all $i$ a.s. Hence, if $M^1:=\langle x_1, X^H\rangle$ and
$T_s:=\inf\{t\ge 0; [ M^1]_t>s\}$, then by Theorem DDS, there exists
a $d$-dimensional Brownian motion $B$ (on an enlarged space) such
that $\overline{Z}=B\circ T$.

Since we already know from part (i) that $Z_t$ has a finite a.s.
limit, it follows that
$$T_{\infty}:=\lim_{t\to\infty} T_t<\infty,\ a.s.\hspace{1cm}\qed$$


\begin{thebibliography}{com}
\bibitem{EHK08}
Engl\"{a}nder, J. , Kyprianou, A. E. and  Harris, S. C.\emph{ Strong
Law of Large Numbers for branching diffusions,} Preprint, 2008.
arXiv:0709.0272

\bibitem{E08}
Engl\"{a}nder, J. \emph{Quenched Law of Large  numbers for Branching
Brownian motion in a random medium}. Ann. Inst. H. Poincar\'e
Probab. Statist.,   44(2008), no. 3,  490-518.

\bibitem{EP99}  Engl\"{a}nder, J. and Pinsky, R. (1999)
\emph{On the construction and support properties of measure-valued
diffusions on $D\subseteq R^{d}$ with spatially dependent
branching}, Ann. Probab. 27(2), 684--730.

\bibitem{K02}  Kallenberg, O. Foundations of modern probability.
Second edition. Probability and its Applications. Springer-Verlag,
New York, 2002.

\bibitem{Ky04} Kyprianou, A. \emph{Asymptotic radial speed of the support of
supercritical branching Brownian motion and super-Brownian motion in
$R\sp d$}. Markov Process. Related Fields 11 (2005), no. 1,
145--156.

\bibitem{P96} Pinsky, R. G.  \emph{Transience, recurrence and local extinction
properties of the support for supercritical finite measure-valued
diffusions.}  Ann. Probab.  24 (1996)(1), 237--267.

\bibitem{T92} Tribe, R.  \emph{The behavior of superprocesses near extinction}.
Ann. Probab. 20 (1992)(1), 286--311.

\end{thebibliography}
\end{document}